\documentclass{amsart}
\usepackage{amssymb}

\usepackage{a4}
\begin{document}
\newtheorem{theorem}{Theorem}[section]
\newtheorem{lemma}[theorem]{Lemma}
\newtheorem{remark}[theorem]{Remark}
\newtheorem{definition}[theorem]{Definition}
\newtheorem{corollary}[theorem]{Corollary}
\newtheorem{example}[theorem]{Example}
\font\pbglie=eufm10
\def\qedbox{\hbox{$\rlap{$\sqcap$}\sqcup$}}
\def\RRR{{\text{\pbglie R}}}
\def\BB{{\mathcal{B}}}
\def\Tr{\operatorname{Tr}}
\def\ffrac#1#2{{\textstyle\frac{#1}{#2}}}
\def\Span{\operatorname{Span}}
\def\Range{\operatorname{Range}}
\def\Rank{\operatorname{Rank}}
\makeatletter
 \renewcommand{\theequation}{%
  \thesection.\alph{equation}}
 \@addtoreset{equation}{section}
 \makeatother
\title[Curvature Homogeneous Pseudo-Riemannian Manifolds]
{Curvature homogeneous spacelike Jordan Osserman pseudo-Riemannian manifolds}
\author{P. Gilkey and S. Nik\v cevi\'c}
\begin{address}{PG: Mathematics Department, University of Oregon,
Eugene Or 97403 USA.\newline Email: {\it gilkey@darkwing.uoregon.edu}}
\end{address}
\begin{address}{SN: Mathematical Institute, Sanu,
Knez Mihailova 35, p.p. 367,
11001 Belgrade,
Yugoslavia.
%stana@sfb288.math.tu-berlin.de
%enikcevi@ubbg.etf.bg.ac.yu
Email: {\it stanan@mi.sanu.ac.yu}}\end{address}
\begin{abstract} Let $s\ge2$. We construct Ricci flat pseudo-Riemannian manifolds of signature
$(2s,s)$ which are not locally homogeneous but whose curvature tensors never the less exhibit
a number of important symmetry properties. They are curvature homogeneous; their curvature tensor is modeled on that of a
local symmetric space. They are spacelike
Jordan Osserman with a Jacobi operator which is nilpotent of order 3; they are not timelike Jordan Osserman. They are
$k$-spacelike higher order Jordan Osserman for $2\le k\le s$; they are
$k$-timelike higher order Jordan Osserman for $s+2\le k\le 2s$, and they are not $k$-timelike higher order Jordan Osserman
for $2\le s\le s+1$.
\end{abstract}
\keywords{Jacobi operator, Osserman manifold, curvature
homogeneous manifold, higher order Osserman manifold, symmetric space.
\newline \phantom{.....}2000 {\it Mathematics Subject Classification.} 53B20.}
%\newline \phantom{.....}\version}
\maketitle

\section{Introduction} Let $\nabla$ be the Levi-Civita connection of a pseudo-Riemannian manifold $(M,g)$ of
signature $(p,q)$ and dimension $m=p+q$. Let
\begin{eqnarray*}
&&R(x,y):=\nabla_x\nabla_y-\nabla_y\nabla_x-\nabla_{[x,y]}\quad\text{and}\\
&&R(x,y,z,w):=g(R(x,y)z,w)
\end{eqnarray*}
be the associated Riemann curvature operator and curvature tensor. Manifolds whose Riemann curvature has a high degree of symmetry
are important in many contexts. Usually this symmetry arises from an underlying symmetry of the metric tensor. One says that $(M,g)$ is {\it
locally homogeneous} if given any two points $P$ and $Q$ of $M$, there exists a local isometry $\psi$ from some neighborhood $U_P$ of $P$ to
some neighborhood $U_Q$ of $Q$ such that $\psi P=Q$.  One says that $(M,g)$ is a {\it local symmetric space} when $\nabla R=0$; local
symmetric spaces are always locally homogeneous. 

In this note, we shall exhibit a family of manifolds whose curvature tensor has a high degree of symmetry
in several different senses, but which are not locally homogeneous. We begin by reviewing some definitions:

\subsection{Curvature homogeneous manifolds} The manifold $(M,g)$ is said to be {\it curvature homogeneous} if given any two
points
$P,Q\in M$, there is a linear isomorphism $\Psi:T_PM\rightarrow T_QM$ so that $\Psi^*g_Q=g_P$
and so that $\Psi^*R_Q=R_P$; see Kowalski, Tricerri,
and Vanhecke \cite{KTV91,KTV92} for further details. If $(M,g)$ is curvature homogeneous, then the curvature tensor looks the
same for every point of $M$.

There is a useful equivalent characterization of curvature homogeneity. Consider the triple
$\mathcal{V}:=(V,g_V,R_V)$ where $g_V$ is a non-degenerate inner product of signature $(p,q)$ on a real vector space $V$
of dimension $m:=p+q$
and
$R_V$ is an algebraic curvature tensor on $V$; i.e. a $4$ tensor which satisfies the usual symmetries of the Riemann curvature tensor:
\begin{eqnarray}
&&R_V(x,y,z,w)=R_V(z,w,x,y)=-R_V(y,x,z,w),\quad\text{and}\label{eqn-1.a}\\
&&R_V(x,y,z,w)+R_V(y,z,x,w)+R_V(z,x,y,w)=0\,.\label{eqn-1.b}
\end{eqnarray}
Equation (\ref{eqn-1.a}) gives $\mathbb{Z}_2$ symmetries; Equation (\ref{eqn-1.b}) is the {\it first Bianchi identity}. We say that
$\mathcal{V}$ is a {\it model space} for
$(M,g,R)$ if given any point $P\in M$, there exists a linear isomorphism
$\Psi:T_PM\rightarrow V$ so $\Psi^*g_V=g_P$ and so $\Psi^*R_V=R_P$; $(M,g)$ is curvature homogeneous if
and only if there exists a model space for $(M,g,R)$. 

\subsection{The Jacobi operator} If $x$ is a tangent
vector at a point $P$ of $M$, then the {\it Jacobi operator} $J(x):y\rightarrow R(y,x)x$ is a self-adjoint endomorphism of the tangent space
$T_PM$. We say that $(M,g)$ is {\it spacelike (\rm resp. \it timelike) Jordan Osserman} if the Jordan normal form of
$J$ is constant on the bundle of unit spacelike (resp. timelike) tangent vectors.

In the Riemannian setting ($p=0$), work of Chi
\cite{refChi} and Nikolayevsky
\cite{refNik,refNik2} shows that if $m\ne8,16$, then any spacelike Jordan Osserman manifold is a $2$ point homogeneous
space; this settles in the affirmative for these dimensions a question raised by Osserman
\cite{refOss}. In the Lorentzian setting ($p=1$), any spacelike or timelike Jordan Osserman manifold necessarily has
constant sectional curvature \cite{refBBG,refGKV}. In the higher signature setting ($p>1,q>1$) the situation is far from clear; we refer to
\cite{refBBGZ,refBBZ,refBV03,refBCGHV,refFeGi03,refGVV,refGiIv} for some partial results.

\subsection{The higher order Jacobi operator} Stanilov and Videv \cite{StVi92} constructed a higher order Jacobi operator. Let $\{e_1,...,e_r\}$
be an orthonormal basis for a spacelike (resp. timelike) $r$-plane $\pi$ in the tangent bundle. The {\it higher order Jacobi operator}
$J(\pi):=J(e_1)+...+J(e_r)$  does not depend on the particular orthonormal basis chosen. One says that $(M,g)$ is {\it $k$-spacelike (\rm
resp. \it $k$-timelike) higher order Jordan Osserman} if the Jordan normal form of $J(\cdot)$ is constant on the Grassmannian
of unoriented spacelike (resp. timelike) $k$-planes. As setting
$k=1$ recovers the previous setting, we shall assume $k\ge2$. The
$k$-spacelike higher order Jordan Osserman manifolds have been classified in the Riemannian setting \cite{refGi02} and in the Lorentzian
setting
\cite{GiSt02} but little is known in the higher signature setting apart from some examples given in \cite{GIZ02}.

\subsection{Curvature homogeneous manifolds which are not locally homogeneous}
It is clear that locally homogeneous manifolds are curvature homogeneous. The somewhat surprising fact is
that the converse fails -- there are curvature homogeneous manifolds which are {\bf not} locally homogeneous. For a
further discussion, we refer to \cite{BK96,BV98,K98,T97,T88,Va91} for Riemannian manifolds,
to \cite{Bu97, BM00,BV97,CLPT90,KM96} for Lorentzian manifolds, and to
\cite{KO99,O96} for affine manifolds. 

In the higher signature setting, there are relatively few examples known of curvature
homogeneous manifolds which are not homogeneous.
The case of signature (2,2) has been studied extensively \cite{BBR00,BCR01,D02,GKV02}; we also refer to \cite{Ha86} for results concerning
isoperimetric hypersurfaces. More generally, let
$p\ge2$. It is known that \cite{refDG,GIZ02,GIZ03} there are pseudo-Riemannian manifolds of neutral signature
$(p,p)$ which are curvature homogeneous but which are not locally homogeneous. These manifolds are
spacelike and timelike Jordan Osserman.
Thus the curvature tensors of these
manifolds exhibit a high degree of symmetry. The Jacobi operator of these manifolds is nilpotent of
order 2.

\subsection{Manifolds of signature $(2s,s)$} Let $s\ge2$. In this paper, we shall extend previous work \cite{GN03} to create {\bf new} examples of
pseudo-Riemannian manifolds of signature $(2s,s)$ whose Riemann curvature tensor also has a high degree of symmetry. The manifolds in this family
are all curvature homogeneous with curvature tensor modeled on that of a symmetric space. Generic members of the family are not
locally homogeneous. They are spacelike Jordan Osserman but not timelike Jordan Osserman. They are $k$-spacelike higher order Jordan Osserman
for
$2\le k\le s$; they are
$k$-timelike higher order Jordan Osserman if and only if
$s+2\le k\le 2s$. Their Jacobi operators are nilpotent of order $3$.

Fix $s\ge2$. We define the pseudo-Riemannian manifolds that we shall be studying and their associated curvature model as follows:

\begin{definition}\label{defn-1.1}\rm
 Let $\vec u:=(u_1,...,u_s)$, $\vec t:=(t_1,...,t_s)$, and $\vec v:=(v_1,...,v_s)$ give
coordinates $(\vec u,\vec t,\vec v)$ on $\mathbb{R}^{3s}$ for $s\ge2$. Let 
$$F(\vec u):=f_1(u_1)+...+f_s(u_s)$$ 
be a smooth function
on an open subset $\mathcal{O}\subset\mathbb{R}^s$. Let 
$$|u|^2:=\textstyle\sum_{1\le i\le s}u_i^2\quad\text{and}\quad u\cdot t:=\sum_{1\le i\le s}u_it_i\,.$$
Define a
pseudo-Riemannian metric
$g_F$ of signature
$(2s,s)$ on $M_F:=\mathcal{O}\times\mathbb{R}^{2s}$ whose non-zero components are given by:
$$
\begin{array}{l}
g_F(\partial_i^u,\partial_i^u)=-2F(\vec u)-2u\cdot t,\\
 g_F(\partial_i^u,\partial_i^v)=g_F(\partial_i^v,\partial_i^u)=1,\\
 g_F(\partial_i^t,\partial_i^t)=-1\,.
\end{array}$$
Set $F_{/j}:=\partial_j^uF=\partial_j^uf_j$, $F_{/ij}:=\partial_j^uF_{/i}$, etc. Note that
$F_{/ij}=0$ for $i\ne j$. Define
\begin{equation}\label{eqn-1.zz}
\alpha_F:=\textstyle\sum_{1\le i\le s}\{F_{/iii}+4u_i\}^2\,.
\end{equation}
\end{definition}

\begin{definition}\label{defn-1.2}
\rm  Let $\{U_1,...,U_s,T_1,...,T_s,V_1,...,V_s\}$ be a basis for $\mathbb{R}^{3s}$ where
$s\ge2$. Let $\mathcal{V}_{3s}:=(\mathbb{R}^{3s},g_{3s},R_{3s})$
where the non-zero entries of the metric $g_{3s}$ and of the algebraic curvature tensor
$R_{3s}$, up to the $\mathbb{Z}_2$ symmetries of Equation (\ref{eqn-1.a}), are
\begin{equation}\label{eqn-1.c}
\begin{array}{l}
g_{3s}(U_i,V_i)=g_{3s}(V_i,U_i)=1,\quad
g_{3s}(T_i,T_i)=-1,\quad\text{and}\\
R_{3s}(U_i,U_j,U_j,T_i)=1\quad\text{for}\quad i\ne j\,.
\end{array}\end{equation}
Set $Z_i^\pm:=U_i\pm\ffrac12V_i$. Then $\Span\{Z_i^+\}$ is a maximal
spacelike subspace of $\mathbb{R}^{3s}$ and $\Span\{T_i,Z_i^-\}$ is the complementary maximal timelike subspace. Thus $\mathbb{R}^{3s}$ has
signature $(2s,s)$. A basis $\BB=\{\tilde U_1,...,\tilde U_s,\tilde T_1,...,\tilde T_s,\tilde V_1,...,\tilde V_s\}$ for $\mathbb{R}^{3s}$ is
said to be {\it normalized} if the relations given above in Display (\ref{eqn-1.c}) hold for $\BB$. 
\end{definition}

\begin{theorem}\label{thm-1.3} Let $s\ge2$. The manifold $(M_F,g_F)$ is a pseudo-Riemannian
manifold of signature $(2s,s)$ which is Ricci flat. We have:
\begin{enumerate}
\item $(M_F,g_F)$ is curvature homogeneous with model $\mathcal{V}_{3s}$.
\item $(M_F,g_F)$ is spacelike Jordan Osserman but not timelike Jordan Osserman.
\item $(M_F,g_F)$ is $k$-spacelike higher order Jordan Osserman for $2\le k\le s$; $(M_F,g_F)$ is $k$-timelike higher order Jordan Osserman if
and only if\newline
$s+2\le k\le 2s$.
\item If there exists a local isometry $\psi$ of $(M_F,g_F)$ with $\psi(P)=Q$, then one has $\alpha_F(P)=\alpha_F(Q)$. Thus
$(M_F,g_F)$ is not locally homogeneous for generic $F$.
\end{enumerate}
\end{theorem}

\subsection{Outline of the paper} In Section \ref{Sect2}, we establish Assertion (1) of Theorem \ref{thm-1.3}
by determining $R$ and $\nabla R$ for $(M_F,g_F)$. It will follow that $\alpha_F$ vanishes identically if and only if $(M_F,g_F)$ is a
local symmetric space. By choosing $F$ suitably, one sees then that $\mathcal{V}_{3s}$ is the model
for the curvature tensor of a Ricci flat local symmetric space. In Sections
\ref{Sect3} and \ref{Sect-4}, we prove Assertions (2) and (3) of Theorem \ref{thm-1.3} by establishing the corresponding results for
the model space
$\mathcal{V}_{3s}$. In Section \ref{Sect5}, we
complete the proof of Theorem
\ref{thm-1.3} by constructing additional natural structures on the manifold
$(M_F,g_F)$ that show $\alpha_F$ is preserved by local isometries.

\section{The curvature tensor of the manifold $(M_F,g_F)$}\label{Sect2}

We begin our study of the manifold $(M_F,g_F)$ by showing:

\begin{lemma}\label{lem-2.1} Let $R_F$ and $\nabla R_F$ be the curvature tensor and the covariant derivative curvature tensor of
the pseudo-Riemannian manifold $(M_F,g_F)$ defined above. Then: 
\begin{enumerate}
\item The non-zero entries in
$R_F$ and
$\nabla R_F$ are, up to the usual $\mathbb{Z}_2$ symmetries,
\begin{enumerate}
\item $R_F(\partial_i^u,\partial_j^u,\partial_j^u,\partial_i^u)
   =F_{/ii}+F_{/jj}+|u|^2$.
\item $R_F(\partial_i^u,\partial_j^u,\partial_j^u,\partial_i^t)=1$.
\item $\nabla R_F(\partial_i^u,\partial_j^u,\partial_j^u,\partial_i^u;\partial_i^u)
  =F_{/iii}+4u_i$.
\end{enumerate}
\item If $\{z_1,...,z_6\}$ are tangent vectors, then $R_F(z_1,z_2)R_F(z_3,z_4)R_F(z_5,z_6)=0$.
\item If $z$ is a tangent vector, then $J_F(z)^3=0$.
\item The manifold $(M_F,g_F)$ is Ricci flat.
\end{enumerate}\end{lemma}

\begin{proof} Let $i\ne j$. The non-zero Christoffel symbols of the second kind are given by:
$$\begin{array}{ll}
g_F(\nabla_{\partial_i^u}\partial_i^u,\partial_i^u)=-F_{/i}-t_i,& \\
g_F(\nabla_{\partial_i^u}\partial_i^u,\partial_j^u)=F_{/j}+t_j,&
g_F(\nabla_{\partial_i^u}\partial_j^u,\partial_i^u)=
g_F(\nabla_{\partial_j^u}\partial_i^u,\partial_i^u)=-F_{/j}-t_j,\vphantom{\vrule height 11pt}\\
g_F(\nabla_{\partial_i^u}\partial_i^u,\partial_i^t)=u_i,&
g_F(\nabla_{\partial_i^u}\partial_i^t,\partial_i^u)=g_F(\nabla_{\partial_i^t}\partial_i^u,\partial_i^u)=-u_i,\\
g_F(\nabla_{\partial_i^u}\partial_i^u,\partial_j^t)=u_j,&
g_F(\nabla_{\partial_i^u}\partial_j^t,\partial_i^u)=
g_F(\nabla_{\partial_j^t}\partial_i^u,\partial_i^u)=-u_j\,.\vphantom{\vrule height 11pt}
\end{array}$$
We may then raise indices to see the non-zero covariant derivatives are given by:
\begin{eqnarray*}
&&\nabla_{\partial_i^u}\partial_i^u=-(F_{/i}+t_i)\partial_i^v
  +\textstyle\sum_{k\ne i}(F_{/k}+t_k)\partial_k^v-\textstyle\sum_{1\le k\le s}u_k\partial_k^t,\\
&&\nabla_{\partial_i^u}\partial_j^u=-(F_{/j}+t_j)\partial_i^v-(F_{/i}+t_i)\partial_j^v,\\
&&\nabla_{\partial_i^u}\partial_i^t=\nabla_{\partial_i^t}\partial_i^u=-u_i\partial_i^v,\quad\text{and}\\
&&\nabla_{\partial_i^u}\partial_j^t=\nabla_{\partial_j^t}\partial_i^u=-u_j\partial_i^v\,.
\end{eqnarray*}

We have $\nabla\partial_i^v=0$. Thus if at least one $z_\mu\in\{\partial_i^v\}$, then
$R_F(z_1,z_2,z_3,z_4)=0$. Similarly, if at least two of the $z_\mu$ belong to $\{\partial_i^t\}$, then $R_F(z_1,z_2,z_3,z_4)=0$. Finally, as
$F_{/ij}=0$ for $i\ne j$,
$R_F(\partial_i^u,\partial_j^u,\partial_k^u,\star)=0$ if the indices $\{i,j,k\}$ are distinct.

The interaction term $-\sum_{1\le k\le s}u_k\partial_k^t$ in $\nabla_{\partial_i^u}\partial_i^u$ is in many ways the
crucial term. We prove Assertions (1a) and (1b) by computing:
$$
\nabla_{\partial_i^u}\nabla_{\partial_j^u}\partial_j^u=F_{/ii}\partial_i^v-\partial_i^t+|u|^2\partial_i^v\quad\text{and}\quad
\nabla_{\partial_j^u}\nabla_{\partial_i^u}\partial_j^u=-F_{/jj}\partial_i^v\,.
$$

We have similarly that $\nabla R_F(X_1,X_2,X_3,X_4;X_5)=0$ if at least one of the $X_i$ belongs to $\Span\{T_i,V_i\}$. Furthermore, up to
the usual
$\mathbb{Z}_2$ symmetries, the only non-zero component of $\nabla R_F$ is given by:
\begin{eqnarray*}
&&\nabla R_F(\partial_i^u,\partial_j^u,\partial_j^u,\partial_i^u;\partial_i^u)\\
&=&\partial_i^uR_F(\partial_i^u,\partial_j^u,\partial_j^u,\partial_i^u)-2
 R_F(\nabla_{\partial_i^u}\partial_i^u,\partial_j^u,\partial_j^u,\partial_i^u)
 -2R_F(\partial_i^u,\nabla_{\partial_i^u}\partial_j^u,\partial_j^u,\partial_i^u)\\
&=&F_{/iii}+2u_i+2R_F(\textstyle\sum_{1\le k\le s}u_k\partial_k^t,\partial_j^u,\partial_j^u,\partial_i^u)+0
=F_{/iii}+4u_i\,.
\end{eqnarray*}
This establishes Assertion (1c). 

Assertions (2) and (3) follow from Assertion (1). Since $J_F(z)^3=0$, $0$ is the only eigenvalue of $J_F(z)$. Thus
$\rho_F(z,z):=\Tr(J_F(z))=0$ and $(M_F,g_F)$ is Ricci flat.
\end{proof}

\begin{proof}[Proof of Theorem \ref{thm-1.3} (1)] Fix $P\in M_F$. Let constants $\varepsilon_i$ and $\varrho_i$ be given. 
We define a new basis
for
$T_PM$ by setting:
$$U_i:=\partial_i^u+\varepsilon_{i}\partial_i^t+\varrho_{i}\partial_i^v,\quad
 T_i:=\partial_i^t+\varepsilon_{i}\partial_i^v,\quad\text{and}\quad V_i:=\partial_i^v\,.$$
Let $i\ne j$. Since $g_F(U_i,T_i)=\varepsilon_{i}-\varepsilon_{i}=0$, the possibly non-zero entries of $g_F$ and $R_F$ are, up to the usual
$\mathbb{Z}_2$ symmetries, given by
\begin{eqnarray*}
&&g_F(U_i,U_i)=g_F(\partial_i^u,\partial_i^u)-\varepsilon_{i}^2+2\varrho_{i},\\
&&g_F(T_i,T_i)=-1,\quad g_F(U_i,V_i)=1,\\
&&R_F(U_i,U_j,U_j,T_i)=1,\quad\text{and}\\
&&R_F(U_i,U_j,U_j,U_i)=F_{/ii}+F_{/jj}+|u|^2+2\varepsilon_i+2\varepsilon_j\,.
\end{eqnarray*}
We set
$$\varepsilon_{i}:=-\ffrac12F_{/ii}-\ffrac14|u|^2
 \quad\text{and}\quad\varrho_i:=\ffrac12\{\varepsilon_{i}^2-g_F(\partial_i^u,\partial_i^u)\}\,.
$$
This ensures that $g_F(U_i,U_i)=0$ and
$R_F(U_i,U_j,U_j,U_i)=0$ and establishes the existence of a basis with the normalizations of Definition \ref{defn-1.2}.
\end{proof}

\begin{remark}\label{rmk-2.2}
\rm Note as a useful scholium to the computations performed above that we can express the function $\alpha_F$ of Equation (\ref{eqn-1.zz}) in
the form:
\begin{eqnarray*}
\alpha_F&=&\ffrac14\textstyle\sum_{i,j,k,l,n}\{\nabla R_F(\partial_i^x,\partial_j^x,\partial_k^x,\partial_l^x;\partial_n^x)\}^2
\\&=&\ffrac14\textstyle\sum_{i,j,k,l,n}\{\nabla R_F(U_i,U_j,U_k,U_l;U_n)\}^2\,.
\end{eqnarray*}
Thus $\nabla R_F$ vanishes if and only if $\alpha_F=0$.
Were one to take $f_i:=-\frac16u_i^3$, then $\alpha_F$ would vanish identically. Consequently, there exist local symmetric spaces
in the family we are considering.\end{remark}

\section{The Jacobi operator of $\mathcal{V}_{3s}$}\label{Sect3} 

In light of Theorem \ref{thm-1.3} (1), one sees that Assertion (2) of Theorem \ref{thm-1.3} will follow
from the corresponding assertions for the model space $\mathcal{V}_{3s}$:

\begin{lemma}\label{lem-3.1} 
Let $J_{3s}$ be the Jacobi operator of $R_{3s}$ for $s\ge2$.
\begin{enumerate}
\item If $g_{3s}(X,X)>0$,
$\Rank\{J_{3s}(X)\}=2(s-1)$ and
$\Rank\{J_{3s}(X)^2\}=s-1$.
\item If $X$ is any element of $\mathbb{R}^{3s}$,
$J_{3s}(X)^3=0$. 
\item The model space $\mathcal{V}_{3s}$ is spacelike Jordan Osserman.
\item The model space $\mathcal{V}_{3s}$ is not timelike Jordan Osserman.
\end{enumerate}\end{lemma}

\begin{proof} There is an additional useful symmetry which plays a crucial role. Let 
$$O(s):=\{\xi=(\xi_{ij}):\textstyle\sum_{1\le i\le s}\xi_{ij}\xi_{ik}=\delta_{jk}\}\subset\mathbb{M}_s(\mathbb{R})$$
 be the standard orthogonal group of $s\times s$ real matrices. Define a diagonal action of
$O(s)$ on $\mathbb{R}^{3s}$ which preserves the structures $g_{3s}$ and $R_{3s}$ by setting:
\begin{equation}\label{eqn-3.a}
\xi:U_i\rightarrow\textstyle\sum_j\xi_{ij}U_j,\quad
 \xi:T_i\rightarrow\textstyle\sum_j\xi_{ij}T_j,\quad\text{and}\quad
 \xi:V_i\rightarrow\textstyle\sum_j\xi_{ij}V_j\,.
\end{equation}

Let $X$ be a spacelike vector. By applying a symmetry of the form described in Equation
(\ref{eqn-3.a}), we may assume that
$$X=a_1U_1+\textstyle\sum_{1\le i\le s}\{b_iT_i+c_iV_i\}\quad\text{where}\quad2a_1c_1-\textstyle\sum_{1\le i\le s}b_i^2>0\,.$$
Thus $a_1\ne0$. Let $i\ge 2$. There exist real numbers $\varepsilon_{ik}\in\mathbb{R}$, where $\varepsilon_{ik}=\varepsilon_{ik}(a,b,c)$ plays
no role in the subsequent development, so that
$$\begin{array}{lll}
J_{3s}(X):X\rightarrow 0,&
J_{3s}(X):T_1\rightarrow 0,&
J_{3s}(X):V_1\rightarrow 0,\\
J_{3s}(X):U_i\rightarrow -a_1^2T_i-\textstyle\sum_{1\le k\le s}\varepsilon_{ik}V_k,&
J_{3s}(X):T_i\rightarrow a_1^2V_i,&
J_{3s}(X):V_i\rightarrow 0\,.
\end{array}$$
This establishes Assertion (1). Assertion (2) is immediate from the definition and
Assertion (3) follows from Assertions (1) and (2).
To establish Assertion (4), we note that
$Z_1^-:=U_1-\frac12V_1$ is a unit timelike vector with $J_{3s}(Z_1^-)\ne0$.
On the other hand,
$T_1$ is also a unit timelike vector with $J_{3s}(T_1)=0$. Thus the Jordan
normal form of $J_{3s}$ is not constant on the pseudo-sphere of unit timelike vectors. \end{proof}

\section{The higher order Jacobi Operator of $\mathcal{V}_{3s}$}\label{Sect-4}

We establish Assertion (3) of Theorem \ref{thm-1.3} by proving:

\begin{lemma}\label{lem-4.1}
Let $J_{3s}$ be the Jacobi operator of $R_{3s}$ for $s\ge2$.
\begin{enumerate}
\item If $\pi$ is a spacelike $k$-plane for $2\le k\le s$, then $\operatorname{Rank}\{J_{3s}(\pi)\}=2s$,
 $\operatorname{Rank}\{J_{3s}(\pi)^2\}=s$, and $J_{3s}(\pi)^3=0$.
\item The model space $\mathcal{V}_{3s}$ is $k$-spacelike higher order Jordan Osserman for\newline $2\le k\le s$.
\item If $\pi$ is a timelike $k$-plane for $s+2\le k\le 2s$, then $\operatorname{Rank}\{J_{3s}(\pi)\}=2s$,
 $\operatorname{Rank}\{J_{3s}(\pi)^2\}=s$, and $J_{3s}(\pi)^3=0$.
\item The model space $\mathcal{V}_{3s}$ is $k$-timelike higher order Jordan Osserman for\newline
$s+2\le k\le 2s$.
\item The model space $\mathcal{V}_{3s}$ is not $k$-timelike higher order Jordan Osserman if\newline $2\le k\le s+1$.
\end{enumerate}
\end{lemma}

\begin{proof} Let $B_U:=\mathbb{R}^{3s}/\Span\{T_i,V_i\}$ and let $\sigma_U$ be the natural projection from $\mathbb{R}^{3s}$ to
$B_U$. Fix a normalized basis $\BB=\{U_i,T_i,V_i\}$ for $\mathbb{R}^{3s}$. Define a positive definite inner product $g_{U,\BB}$ on
$B_U=\Span\{\sigma_U(U_i)\}$ so
$$g_{U,\BB}(\sigma_UU_i,\sigma_UU_j)=\delta_{ij}\,.$$
We will show presently in Lemma \ref{lem-7.1} that $B_U$ and $g_{U,\BB}$ are independent of the particular normalized basis which was chosen,
but this plays no role at present. If $\pi$ is a linear subspace of $\mathbb{R}^{3s}$, set
$$
\tilde g_\pi:=\sigma_U^*g_{U,\BB}|_\pi\quad\text{and}\quad\ell(\pi):=\Rank\{\tilde g_\pi\}\,.
$$

Let $\pi$ be a spacelike $k$-plane where $k\ge2$. Since every non-zero vector of $\pi$ is spacelike, $\pi\cap\ker(\sigma_U)=\{0\}$. Thus
$\ell(\pi)=k$. Let indices $\alpha,\beta$ range from $1$ thru $k$; let the index $\mu$ range from $k+1$ thru $s$. By diagonalizing the inner
product
$\tilde g_\pi$ with respect to the positive definite inner product
$g_{3s}|_\pi$, we can choose an orthonormal basis
$\{X_\alpha\}$ for
$\pi$ and {\it positive constants} $a_\alpha$ so that
$$
\tilde g_\pi(X_\alpha,X_\beta)=a_\alpha a_\beta\delta_{\alpha\beta}\quad\text{and}\quad
g_{3s}(X_\alpha,X_\beta)=\delta_{\alpha\beta}\,.
$$
By applying a symmetry of the form described in Equation (\ref{eqn-3.a}), we may suppose
$$
X_\alpha=a_\alpha U_\alpha+\textstyle\sum_j\{b_{\alpha j}T_j+c_{\alpha j}V_j\}\quad
\text{and}\quad g_{3s}(X_\alpha,X_\beta)=\delta_{\alpha\beta}
$$
for suitably chosen constants $b_{\alpha j}$ and $c_{\alpha j}$. As $J_{3s}(\pi)=\sum_\alpha J_{3s}(X_\alpha)$, there exist
constants $\varepsilon_{ij}=\varepsilon_{ij}(a,b,c)\in\mathbb{R}$ so
$$\begin{array}{ll}
J_{3s}(\pi):U_\beta\rightarrow-\textstyle\sum_{\alpha\ne\beta}a_\alpha^2T_\beta+\sum_j\varepsilon_{\beta j}V_j,&
J_{3s}(\pi):T_\beta\rightarrow\textstyle\sum_{\alpha\ne\beta}a_\alpha^2V_\beta,\\
J_{3s}(\pi):U_\nu\rightarrow-\textstyle\sum_\alpha a_\alpha^2T_\nu+\sum_j\varepsilon_{\nu j}V_j,&
J_{3s}(\pi):T_\nu\rightarrow\textstyle\sum_\alpha a_\alpha^2V_\nu,\vphantom{\vrule height 11pt}\\
J_{3s}(\pi):V_\beta\rightarrow0,&
J_{3s}(\pi):V_\nu\rightarrow0\,.\vphantom{\vrule height 11pt}
\end{array}$$
Since $k\ge2$, one has that $\sum_{\beta\ne\alpha}a_\beta^2\ne0$. Consequently
\begin{eqnarray*}
&&\Range J_{3s}(\pi)=\Span\{T_1,...,T_s,V_1,...,V_s\},\\
&&\Range J_{3s}(\pi)^2=\Span\{V_1,...,V_s\},\quad\text{and}\quad\Range J_{3s}(\pi)^3=\{0\}\,.
\end{eqnarray*}
Assertion (1) now follows. Assertion (2) follows from Assertion (1).

Let $\pi$ be a timelike $k$-plane. We diagonalize $\tilde g_\pi$ with respect to the negative definite quadratic form $g_{3s}|_\pi$ to choose an
orthonormal basis $\{X_\alpha\}$ for $\pi$ so that 
$$\tilde g_\pi(X_\alpha,X_\beta)=a_{\alpha}a_\beta\delta_{\alpha\beta}\quad\text{and}\quad
  g_{3s}(X_\alpha,X_\beta)=-\delta_{\alpha\beta}\,.$$
We have $\ell(\pi)$ is the number of times that $a_\alpha\ne0$. Again, by applying an
appropriate symmetry $\xi\in O(s)$ as described in Equation (\ref{eqn-3.a}), we can assume without loss of generality
$$X_\alpha=a_\alpha U_\alpha+\textstyle\sum_{1\le i\le s}\{b_{\alpha i}T_i+c_{\alpha i}V_i\}\quad\text{for}\quad1\le\alpha\le k\,.$$
The calculations performed above show that if $\ell\ge2$, then
\begin{equation}\label{eqn-4.a}
\Rank\{J_{3s}(\pi)\}=2s,\quad\Rank\{J_{3s}(\pi)^2\}=s,\quad\text{and}\quad J_{3s}(\pi)^3=0\,.
\end{equation}

Since $\ker(\sigma_U)=\Span\{T_i,V_i\}$, any timelike subspace of $\ker(\sigma_U)$ has dimension at most $s$. Since $\pi$ is
timelike, $\dim\{\pi\cap\ker(\sigma_U)\}\le s$ and hence
$$
\ell=\dim\{\sigma_U(\pi)\}=\dim\{\pi\}-\dim\{\pi\cap\ker(\sigma_U)\}\ge k-s\,.
$$
Thus if $k\ge s+2$, then $\ell\ge2$. Assertion (3) now follows from Equation (\ref{eqn-4.a}); Assertion (4) follows from Assertion (3).

To prove the final assertion, we must give examples of timelike $k$-planes whose Jacobi operators have different Jordan normal forms. The
calculations performed above show that:
$$\Rank\{J_{3s}(\pi)\}=\left\{\begin{array}{lll}
0&\text{if}&\ell(\pi)=0,\\
s-1&\text{if}&\ell(\pi)=1,\\
s&\text{if}&\ell(\pi)\ge2\,.
\end{array}\right.$$ 
If
$2\le k\le s+1$, set
\begin{eqnarray*}
&&\pi_1:=\left\{\begin{array}{ll}
\Span\{T_1,...,T_k\}&\text{if }k\le s,\\
\Span\{T_1,...,T_s,Z_1^-\}\hphantom{.......,a}&\text{if }k=s+1,
\end{array}\right.\\
&&\pi_2:=\left\{\begin{array}{ll}
\Span\{T_1,...,T_{k-1},Z_1^-\}&\text{if }k\le s,\\
\Span\{T_1,...,T_{s-1},Z_1^-,Z_2^-\}&\text{if }k=s+1\,.
\end{array}\right.\end{eqnarray*}
Then $\pi_1$ and $\pi_2$ are timelike $k$-planes with
\begin{eqnarray*}
&&\Rank\{J_{3s}(\pi_1)\}=\left\{\begin{array}{lll}
0&\text{if }k\le s&\text{as }\ell(\pi)=0,\\
s-1&\text{if }k=s+1&\text{as }\ell(\pi)=1,\end{array}\right.\\
&\ne&\Rank\{J_{3s}(\pi_2)\}=\left\{\begin{array}{lll}
s-1&\text{if }k\le s&\text{as }\ell(\pi)=1,\\
s&\text{if }k=s+1&\text{as }\ell(\pi)=2.\end{array}\right.
\end{eqnarray*}
Consequently $\mathcal{V}_{3s}$ is not $k$-timelike higher order Jordan Osserman.
\end{proof}

\section{Invariants of the manifold $(M_F,g_F)$}\label{Sect5}

It is clear from the definition that $||R||_{g_F}=0$ and $||\nabla R||_{g_F}=0$. Thus to prove the final assertion of Theorem
\ref{thm-1.3}, we must introduce some additional structures and show that they are invariantly defined. We work on the model
space $\mathcal{V}_{3s}$ and suppress the index $s$ in the interests of notational simplicity. We now return to structures introduced earlier in
the proof of Lemma \ref{lem-4.1} and show these structures are intrinsic -- i.e. they are independent of the particular normalized basis which
was chosen. Consider the following subspaces of $\mathbb{R}^{3s}$:
\begin{eqnarray*}
&&A_V:=\{W\in\mathbb{R}^{3s}:R(W_1,W_2,W_3,W)=0\ \forall\ W_1,W_2,W_3\in\mathbb{R}^{3s}\},\\
&&A_{T,V}:=A_V^\perp=\{W\in\mathbb{R}^{3s}:g(W,W_1)=0\ \forall\ W_1\in A_V\}\,.
\end{eqnarray*}
Let $\sigma_{U,T}$, $\sigma_T$, and $\sigma_U$ be the natural projections to the quotient spaces
$$
B_{U,T}:=\mathbb{R}^{3s}/A_V,\quad
B_T:=A_{T,V}/A_V,\quad\text{and}\quad
B_U:=\mathbb{R}^{3s}/A_{T,V}\,.
$$

The spaces given above are defined invariantly; they do not depend on the choice of basis. On the other hand, if $\BB:=\{
U_i, T_i, V_i\}$ is {\bf any} basis for
$\mathbb{R}^{3s}$ which satisfies the normalizations given in Definition \ref{defn-1.2}, then one may express:
$$\begin{array}{ll}
A_V=\Span\{ V_i\},&A_{T,V}=\Span\{ T_i, V_i\},\\
B_{U,T}=\Span\{\sigma_{U,T} U_i,\sigma_{U,T} T_i\},&B_T=\Span\{\sigma_T T_i\},
 \vphantom{\vrule height 11pt}\\
B_U=\Span\{\sigma_U U_i\}\,.\vphantom{\vrule height 11pt}
\end{array}$$

The metric $g_{3s}$ descends to a negative definite inner product $g_T$ on $B_T\subset B_{U,T}$;
$\{\sigma_T( T_i)\}$ is an orthonormal basis for $B_T$. Note that $g_T$ is {\bf not} defined
on all of $B_{U,T}$ but only on the subspace $B_T$. 
 Let $g_{U,\BB}(\sigma_U U_i,\sigma_U U_j)=\delta_{ij}$ define a positive
definite metric $g_{U,\BB}$ on
$B_U$ which a priori depends on the basis $\BB$. 

\begin{lemma}\label{lem-7.1}We have $g_{U,\BB}=g_{U,\tilde\BB}$ for any two normalized bases $\BB$ and $\tilde\BB$ of
$\mathbb{R}^{3s}$.
\end{lemma}

\begin{proof} The tensor $R$ descends to a tensor $R_{U,T}$ on $B_{U,T}$ so that $\sigma_{U,T}^*R_{U,T}=R$. The basis dependent action
of the orthogonal group $O(s)$ on $T_PM$ described in Equation (\ref{eqn-3.a}) induces basis dependent actions on the subspaces
$A_V$ and $A_{T,V}$
and on the quotient spaces $B_U$, $B_T$, and $B_{U,T}$ described above. This action preserves the metric $g_T$, the metric $g_{U,\BB}$, and the
tensor $R_{U,T}$.

Let $\BB$ and $\tilde\BB$ be normalized bases for $\mathbb{R}^3$. Then $\{\sigma_T(\tilde T_i)\}$ and $\{\sigma_T(T_i)\}$ are orthonormal
bases for
$B_T$. By replacing
$\tilde\BB$ by
$\xi\tilde\BB$ if necessary, where $\xi$ is a suitably chosen element of $O(s)$, we can assume without loss of generality
$\sigma_{U,T}\tilde T_i=\sigma_{U,T}T_i$ for all
$i$.  Let $u_i=\sigma_{U,T}U_i$, $\tilde u_i:=\sigma_{U,T}\tilde U_i$, and $t_i:=\sigma_{U,T}T_i=\sigma_{U,T}\tilde T_i$. 
Expand $$\tilde u_j=\textstyle\sum_{1\le k\le s}\{a_{jk}u_k+b_{jk}t_k\}\,.$$ We shall prove the Lemma by showing that $a_{jk}=\delta_{jk}$.

Let $j\ne k$. We use the defining relations to see
\begin{eqnarray*}
&&1=R_{U,T}(\tilde u_{j},\tilde u_{k},\tilde u_{k},t_{j})
  =(a_{jj}a_{kk}-a_{jk}a_{kj})a_{kk}\\
&&0=R_{U,T}(\tilde u_{j},\tilde u_{k},\tilde u_{j},t_{j})
 =(a_{jj}a_{kk}-a_{jk}a_{kj})a_{jk}\,.
\end{eqnarray*}
Since $0\ne(a_{jj}a_{kk}-a_{jk}a_{kj})$, we have $a_{jk}=0$ for $j\ne k$; similarly $a_{kj}=0$ for $j\ne k$. Thus $a_{jj}a_{kk}a_{kk}=1$.
Similarly
$a_{jj}a_{kk}a_{jj}=1$. Thus $a_{kk}=1$ so $a_{jk}=\delta_{jk}$.
\end{proof}

\begin{proof}[Proof of Theorem \ref{thm-1.3} (4).] Fix $P\in TM$. Let
$\{\partial_i^u,\partial_i^t,\partial_i^v\}$ be the coordinate frame for $T_PM$. We use the adjusted basis $\{U_i,T_i,V_i\}$ constructed in
Section
\ref{Sect2} to find an isomorphism $\Psi$ which identifies $(T_PM,g_F,R_F)$ with $\mathcal{V}_{3s}$. As
$\nabla R_F(\star,\star,\star,\star;\star)=0$ if any entry belongs to
$A_{T,V}$, there is a tensor $\nabla R_{U}$ on $B_U$ so $\nabla R=\Psi^*\sigma_U^*\nabla R_{U}$. Let $\alpha_F$ be as
defined previously. We use Remark \ref{rmk-2.2} to see
\begin{eqnarray*}
&&\alpha_P=\ffrac14\sum_{i_1,i_2,i_3,i_4,i_5}\nabla
R(\partial_{i_1}^u,\partial_{i_2}^u,\partial_{i_3}^u,\partial_{i_4}^u;\partial_{i_5}^u)^2\\
&=&\ffrac14\sum_{i_1,i_2,i_3,i_4,i_5}\nabla R(U_{i_1},U_{i_2},U_{i_3},U_{i_4};U_{i_5})^2
=\ffrac14||\nabla R_{U}||_{g_U}^2\,.
\end{eqnarray*}
As $||\nabla R_U||_{g_U}^2$ is invariantly defined, $\alpha_P$ is preserved by local isometries. Thus, if $(M_F,g_F)$ is
locally homogeneous, then
$\alpha_P$ is constant. This fails for generic $F$.
\end{proof}

\begin{remark}\label{rmk-7.2}
\rm We can construct additional invariants of the metric by considering the norms of
higher order covariant derivatives of the curvature tensor. Set:
\begin{eqnarray*}
\alpha_F^k:&=&2^{-k-1}||\nabla^{(k)}R||_{g_U}\\
&=&2^{-k-1}\sum_{i_1,i_2,i_3,i_4,j_1,...,j_k}
R(\partial_{i_1}^u,\partial_{i_2}^u,\partial_{i_3}^u,\partial_{i_4}^u;\partial_{j_1}^u,...,\partial_{j_k}^u)^2\,.
\end{eqnarray*}
\end{remark}

\section*{Acknowledgments} Research of P. Gilkey partially supported by the
MPI (Leipzig). Research of S. Nik\v cevi\'c partially supported by the DAAD (Germany) and MM 1646 (Srbija). The authors wish to express their thanks to the Technishe Universit\"at
Berlin where much of the research reported here was conducted. Finally, it is a pleasant task to
thank Professor E. Garc\'{\i}a--R\'{\i}o for helpful discussions.


\begin{thebibliography}{AAA}


\bibitem{refBBG} N. Bla\v zi\'c, N. Bokan and P. Gilkey,
  {\it A Note on Osserman Lorentzian manifolds}, Bull. London Math. Soc.
  {\bf 29} (1997), 227--230.



\bibitem{refBBGZ} N. Bla\v zi\'c, N. Bokan, P. Gilkey and Z. Raki\'c,
   {\it Pseudo-Riemannian Osserman manifolds}, Balkan J. Geom. Appl.
   {\bf 2}, (1997), 1--12.

\bibitem{refBBZ} N. Bla\v zi\'c, N. Bokan, and Z. Raki\'c,
   {\it Osserman pseudo-Riemannian manifolds of signature $(2,2)$}, Aust. Math. Soc.
   {\bf 71}, (2001), 367--395.

\bibitem{BBR00}N. Bla\v zi\'c, N. Bokan, and Z. Raki\'c,
{\it A note on the Osserman conjecture and isotropic covariant derivative of curvature},
Proc. Amer. Math. Soc. {\bf 128} (2000), 245--253.

\bibitem{refBV03} N. Bla\v zi\' c and S. Vukmirovi\' c, {\it Examples of self-dual,
Einstein metrics of $(2,2)$-signature}, math.DG/0206081.

\bibitem{BK96} E. Boeckx, O. Kowalski, and L. Vanhecke, {\bf Riemannian manifolds of
conullity two}, World Scientific Publishing Co., Inc., River Edge, NJ, 1996. ISBN 981-02-2768-X.

\bibitem{BV98} E. Boeckx and L. Vanhecke, {\it Curvature homogeneous unit tangent sphere bundles},
Publ. Math. Debrecen {\bf 53} (1998), 389--413.


\bibitem{BCR01} A. Bonome, P. Castro, E. Garc\'{\i}a--R\'{\i}o,
{\it Generalized Osserman four-dimensional manifolds}, Classical Quantum Gravity {\bf 18} (2001), 4813--4822.

\bibitem{refBCGHV} A. Bonome, R. Castro, E. Garc\'{\i}a--R\'{\i}o, L. M. Hervella, R. V\'{a}zquez-Lorenzo,
{\it Nonsymmetric
 Osserman indefinite K\"{a}hler manifolds}, Proc. Amer. Math. Soc. {\bf 126} (1998), 2763--2769.


\bibitem{Bu97} P. Bueken, {\it On curvature homogeneous three-dimensional Lorentzian manifolds}, 
J. Geom. Phys. {\bf 22} (1997), 349--362.


\bibitem{BM00} P. Bueken and M. Djori\'c, 
{\it Three-dimensional Lorentz metrics and curvature homogeneity of order one},
Ann. Global Anal. Geom. {\bf 18} (2000), 85--103. 

\bibitem{BV97} P. Bueken and L. Vanhecke, 
{\it Examples of curvature homogeneous Lorentz metrics}, 
Classical Quantum Gravity {\bf 14} (1997), L93--L96.

\bibitem{CLPT90} M. Cahen, J. Leroy, M. Parker, F. Tricerri, and L. Vanhecke, 
{\it Lorentz manifolds modelled on a Lorentz symmetric space}, 
J. Geom. Phys. {\bf 7} (1990), 571--581.


\bibitem{refChi} Q.-S. Chi, 
{\it A curvature characterization of certain locally rank-one symmetric spaces}, 
J. Differential Geom. {\bf 28} (1988), 187--202.

\bibitem{D02} A. Derdzinski, 
{\it Curvature homogeneous indefinite Einstein metrics in dimension four: the diagonalizable case}, 
to appear Contemporary Mathematics; 
math.DG/0211248.

\bibitem{refDG} C. Dunn and P. Gilkey, 
{\it Curvature homogeneous pseudo-Riemannian manifolds which are not locally homogeneous}; 
math.DG/0306072.

\bibitem{refFeGi03} B. Fiedler and P. Gilkey, 
{\it Nilpotent Szab\'o, Osserman and Ivanov-Petrova pseudo-Riemannian manifolds}, 
to appear Contemporary Mathematics; math.DG/0211080.

\bibitem{refGKV} E. Garc\'{\i}a--R\'{\i}o, D. Kupeli and M. E. V\'azquez-Abal, 
{\it On a problem of Osserman in Lorentzian geometry},
Differential Geom. Appl. {\bf 7} (1997), 85--100.

\bibitem{GKV02} E. Garc\'{\i}a--R\'{\i}o, D. Kupeli, and R. V\'azquez-Lorenzo, 
{\bf Osserman Manifolds in Semi-Riemannian Geometry}, 
Lecture Notes in Mathematics, 1777. 
Springer-Verlag, Berlin, 2002.
ISBN: 3-540-43144-6.

\bibitem{refGVV} E. Garc\'ia-Ri\'o, M. E. V\' azquez-Abal and R. V\' azquez-Lorenzo, 
{\it Nonsymmetric Osserman pseudo-Riemannian manifolds},
  Proc. Amer. Math. Soc. {\bf 126} (1998), 2771--2778.

\bibitem{refGi02} P. Gilkey, {\bf Geometric properties of natural operators defined by the Riemann
curvature tensor}, World Scientific Publishing Co., Inc., River Edge, NJ, 2001. ISBN: 981-02-4752-4.

\bibitem{refGiIv} P. Gilkey and R. Ivanova, {\it Spacelike Jordan Osserman algebraic curvature tensors in the
higher signature setting}, 
{\bf Differential geometry}, Valencia, 2001, 179--186, World Sci. Publishing, River Edge, NJ, 2002.
ISBN 981-02-4906 (2002).

\bibitem{GIZ02} P. Gilkey, R. Ivanova, and T. Zhang, 
 {\it Higher order Jordan Osserman pseudo-Riemannian manifolds}, 
 Classical Quantum Gravity, {\bf 19} (2002), 4543--4551.

\bibitem{GIZ03} ---
%P. Gilkey, R. Ivanova, and T. Zhang, 
 {\it Szab\'o Osserman IP Pseudo-Riemannian manifolds}, 
 Publ. Math. Debrecen 62 (2003), 387--401.

\bibitem{GN03} P. Gilkey and S. Nik\v cevi\'c, 
{\it Nilpotent Spacelike Jorden Osserman pseudo-Riemannian manifolds}; 
math.DG/0302044.

\bibitem{GiSt02} P. Gilkey and I. Stavrov,
{\it Curvature tensors whose Jacobi or Szab\'o operator is nilpotent on null vectors},
 Bull. London Math. Soc. {\bf 34} (2002), 650--658.

\bibitem{Ha86} J. Hahn,
{\it Homogene Hyperfl\"achen in der pseudoremannschen geometrie},
Bonner Mathematische Schriften 172, Universit\"at Bonn Mathematische Institute, Bonn (1986).

\bibitem{KM96} A. Koutras and C. McIntosh, 
 {\it A metric with no symmetries or invariants}, 
 Classical Quantum Gravity {\bf 13} (1996), L47--L49.

\bibitem{K98} O. Kowalski, 
 {\it On curvature homogeneous spaces}, 
 {\bf Proceedings of the workshop on recent topics in differential geometry},
Santiago de Compostela, Spain, July 16--19, 1997. Santiago de Compostela: Universidade de Santiago de
Compostela. Publ. Dep. Geom. Topolog\'a, Univ. Santiago Compostela. 89, 193--205 (1998).

\bibitem{KO99} O. Kowalski, B. Opozda, and Z. Vl\v sek, 
 {\it Curvature homogeneity of affine connections on two-dimensional manifolds}, 
 Colloq. Math. {\bf 81} (1999), 123--139.

\bibitem{KTV91} O. Kowalski, F. Tricerri, and L. Vanhecke, 
{\it New examples of non-homogeneous Riemannian manifolds whose curvature tensor is that of a Riemannian symmetric space}, 
C. R. Acad. Sci. Paris S\'er. I Math. {\bf 311} (1990), 355--360.

\bibitem{KTV92} ---,
{\it Curvature homogeneous Riemannian manifolds},
 J. Math. Pures Appl. {\bf 71} (1992), 471--501.

\bibitem{refNik} Y. Nikolayevsky, 
{\it Two theorems on Osserman manifolds},
Differential Geom. Appl. {\bf 18} (2003), 239--253.

\bibitem{refNik2} Y. Nikolayevsky, Osserman Conjecture in dimension $n \ne 8, 16$,
preprint:\newline
http://arXiv.org/abs/math.DG/0204258.

\bibitem{O96} B. Opozda, 
  {\it On curvature homogeneous and locally homogeneous affine connections}, 
Proc. Amer. Math. Soc. {\bf 124} (1996), 1889--1893.

\bibitem{refOss} R. Osserman, 
  {\it Curvature in the eighties}, 
  Amer. Math. Monthly, {\bf97}, (1990) 731--756.


\bibitem{StVi92} G. Stanilov and V. Videv,
{\it On a generalization of the Jacobi operator in the Riemannian geometry}, Annuaire Univ. Sofia Fac.
Math. Inform. {\bf 86} (1992), 27--34.

\bibitem{St03} I. Stavrov, {\bf Spectral geometry of the Riemann curvature tensor}, Ph. D. Thesis, University of Oregon (2003).

\bibitem{T97} A. Tomassini, 
 {\it Curvature homogeneous metrics on principal fibre bundles},
Ann. Mat. Pura Appl. {\bf 172} (1997), 287--295.

\bibitem{T88} F. Tricerri,
{\it Riemannian manifolds with the same curvature as a homogeneous space, and a conjecture of Gromov}. 
 {\bf Geometry Conference} (Parma, 1988). Riv. Mat. Univ. Parma (4) {\bf 14} (1988), 91--104.\vfill\eject

\bibitem{TzVi99} J.Tzankov and V.Videv, 
{\it A Riemannian pointwise Stanilov manifold of type (n,k)},
Abstracts of the 4th International Conference on Geometry and Applications, Varna, 1999, 62--63.


\bibitem{Va91} L. Vanhecke, {\it Curvature homogeneity and related problems},
{\bf Proceedings of the Workshop on Recent Topics in Differential Geometry} (Puerto de la Cruz, 1990), 103--122, 
Informes, 32, Univ. La Laguna, La Laguna, 1991. 

\end{thebibliography}
\end{document}